\documentclass[12pt]{article}
\usepackage{latexsym, amsfonts}
\usepackage{epsfig}
\newtheorem{thm}{Theorem}[section]
\newtheorem{lem}{Lemma}[section]

\newtheorem{rem}{Remark}[section]
\newtheorem{defn}{Definition}[section]

\def\ex{\mbox{\rm E\,}}
\def\var{\mbox{\rm Var\,}}

\def\I{\mbox{\rm I}}
\def\mse{\mbox{\rm MSE}}

\begin{document}

\title{Estimating a Structural Distribution\\
 Function by Grouping}

\author{Bert van Es and Stamatis Kolios \\
{\normalsize Korteweg-de Vries Institute for Mathematics}\\
{\normalsize University of Amsterdam}\\
{\normalsize The Netherlands}}

\maketitle

\begin{abstract}
By the method of Poissonization we confirm some existing results concerning consistent estimation of the structural
distribution function in the situation of a large number of rare events.
Inconsistency of the so called natural estimator is proved.
The method of grouping in cells of equal size is investigated and its consistency derived.
A bound on the mean squared error is derived.
\\[.3cm]

\end{abstract}

\section{Introduction and results}
\setcounter{equation}{0}

The concept of a structural distribution function originates from linguistics.
Let $M$ denote the size of  the vocabulary of an author and consider a text of
this author that contains $n$ words.
Every choice of a word in the text from the vocabulary can be seen as the realization
of a multinomial random vector. The whole text consists of
a sequence of such choices $X^{(i)}=( X^{(i)}_{1,M},\ldots, X^{(i)}_{M,M}), \hspace{10pt} i=1,2,\ldots,n$
, which are assumed to be independent.
So each $X^{(i)}$ is Multinomial$(1, p_{1,M}, p_{2,M},\ldots,p_{M,M})$ distributed,
where $p_{1,M}, p_{2,M},\ldots,p_{M,M}$
denote the {\it cell probabilities}. In linguistics the vector of those word probabilities is viewed as
a characteristic of the author. More specifically one is interested in estimating the so called
{\em structural distribution function}.

\begin{defn}
The Structural Distribution Function $F_M$ is the empirical distribution
function based on $M$ times the cell probabilities. Hence
\begin{equation}
   F_M(x)=\frac{1}{M} \sum_{j=1}^{M} I_{[Mp_{j,M}\leq x]}.
\end{equation}
\end{defn}
\bigskip
\noindent

We will investigate the estimation problem for the case of a
{\em large number of rare events}, i.e. we assume
\begin{equation}\label{cond}
n, M \to \infty \quad\mbox{and}\quad
n/M \to \lambda,\ \mbox{where}\ 0< \lambda < \infty.
\end{equation}
So in the linguistic context both sizes of the text and the vocabulary
are large, and the text size is proportional to the size of the vocabulary.
Assuming that, under (\ref{cond}), $F_M$ converges
weakly to a distribution function $F$ we want to estimate $F$ at a
fixed positive point $x$.
The problem of estimation of  $p_{1,M}, p_{2,M},\ldots,p_{M,M}$ is thus asymptotically
replaced by estimation of  $F$.

The estimators we consider are based on the {\it cell counts} of the $n$ observations of $X$
, i.e.
\begin{equation}
       \nu_{j,M}=\sum_{i=1}^{n}X_{j,M}^{(i)}, \hspace{3pt}
       j=1,2,\ldots,M.
\end{equation}
Since the cell probabilities can be estimated by the cell frequencies
an obvious estimator of $F$ seems to be the {\it natural estimator}
$\hat F_M $ which is defined as the empirical distribution function based on
$M$ times the {\it cell frequencies} $\nu_{j,M} /n$. Hence

\begin{equation}
   \hat F_M(x)=\frac{1}{M} \sum_{j=1}^{M} I_{[\frac{M}{n}\nu_{j,M}\leq
   x]}.
\end{equation}

The method of Poissonization is based on the following idea. Instead
of considering the cell counts based on $n$ observations of $X$, we
introduce the cell counts $\rho_{j,M}$ based on $N$ observations of $X$,
where $N$ is a Poisson(n) distributed random variable independent of the $X$'s. So
\begin{equation}
       \rho_{j,M}=\sum_{i=1}^{N}X_{j,M}^{(i)}, \hspace{3pt}
       j=1,2,\ldots,M.
\end{equation}
The advantage of Poissonization is  that the $\rho_{j,M}$ are {\em independent}
 Poisson($np_{j,M}$) random
variables, while ($\nu_{1,M},\ldots,\nu_{M,M}$) are
Multinomial($n,p_{1,M},p_{2,M},\ldots,p_{M,M}$) distributed.

The natural estimator based on $\rho_{1,M},\rho_{2,M},\ldots,\rho_{M,M}$, denoted  by $\tilde
F_M(x)$,
is then equal to
\begin{equation}
\tilde F_M(x)=\frac{1}{M} \sum_{j=1}^{M} I_{[\frac{M}{n}\rho_{j,M}\leq
x]}.
\end{equation}

Let $Z_M$ denote a random variable with distribution function $F_M$ and $Z$
a random variable with distribution function $F$. The following theorem establishes
the inconsistency of the natural
estimator. This has already been proved by Klaassen and Mnatsakanov (2000) without using
 Poissonization.
\begin{thm} \label{thm:1.1}
Let (\ref{cond}) hold and let $F_M \stackrel{w}{\to}F$
(or equivalently $Z_M \stackrel{w}{\to}Z$). Then
\begin{equation}\label{lim1}
\hat F_M(x) \stackrel{P}{\to} F_{Y/\lambda}(x),
\end{equation}
where the conditional distribution of $Y$ given $Z=z$ is
Poisson($\lambda z$), for positive $z$,
 and of $Y$ given $Z=0$ is degenerate at zero.
\end{thm}

Inconsistency of $\hat F_M$ also follows from the fact that it is a distribution
function with jumps only at multiples
of $M/n$. Hence, in the limit, it can only have mass at multiples of $1/\lambda$.
However, knowledge of the limit is
useful since based on the exact limit
given by Theorem~\ref{thm:1.1}, Klaassen and Mnatsakanov (2000)
have constructted a consistent estimator of $F$ by Laplace inversion.

\bigskip

The inconsistency of the natural estimator seems to occur since $n$ increases too slowly with
regard to the number of cells $M$. We can reduce that number by replacing the
$M$ cells by $m$ groups and assuming $n/m \to \infty$. We define the  {\it
grouped cell probabilities} $q_{j,M}$ by
\begin{equation}
  q_{j,M}=\sum_{i=k_{j-1}+1}^{k_j}p_{i,M}, \hspace{3pt} j=1,2,\ldots,m
\end{equation}
and the {\it grouped cell frequencies} $\bar \nu_{j,M}$ as
\begin{equation}
  \bar \nu_{j,M}=\sum_{i=k_{j-1}+1}^{k_j}\nu_{i,M}, \hspace{3pt}
j=1,2,\ldots,m
\end{equation}
where the {\it cell limits} $k_j, \hspace{3pt} j=0,1,\ldots,m$, are integers
such that $0=k_0<k_1< \ldots<k_m=M$. We restrict ourselves to the situation
where the $m$ groups are of {\em equal size} $k$, so
$M=km$ and $k_j = j k$.

Let $F_m$ denote the empirical distribution function based on $m$
times the grouped cell probabilities. So
\begin{equation}
F_m(x)={1\over m} \sum_{j=1}^{m}I_{[mq_{j,M} \leq x]}.
\end{equation}
\noindent
Define the estimator $\hat F_m(x)$ based on the grouped cell counts by
\begin{equation}\label{groupest}
\hat F_m(x)={1\over m} \sum_{j=1}^{m} I_{[\frac{m}{n} \bar \nu_{j,M} \leq
x]}.
\end{equation}
The Poissonized version ${\tilde F}_m(x)$, based on the grouped Poisson
counts
\begin{equation}\label{grouppoisson}
\bar \rho_{j,M}= \sum_{i=k_{j-1}+1}^{k_j}\rho_{i,M},\ j=1,\dots,m,
\end{equation}
 is obtained by replacing the $\bar\nu$'s
by $\bar\rho$'s in (\ref{groupest}).
Note that $\bar \rho_{j,M}$ has a Poisson$(n q_{j,M})$ distribution and that the $\bar\rho$'s  are
independent.
Note also that for $m=M$ and hence $k=1$, a situation excluded by condition (\ref{condm})
below, we regain the natural estimator $\hat F_M(x)$.

The following theorem establishes the weak consistency of the
estimator based on the grouped counts.
\begin{thm} \label{thm:1.2}
Let (\ref{cond}) hold. Assume further that
\begin{equation}\label{condm}
{n\over m\log m}\to \infty.
\end{equation}
If $F_m\stackrel{w}{\to} F$  and the distributions induced by the $F_m$ are concentrated
on a fixed bounded set, then
\begin{equation}\label{lim2}
\hat F_m(x) \stackrel{P}{\to} F(x),
\end{equation}
for every  continuity point $x$ of $F.$
\end{thm}

Let us sketch the proofs of the two theorems.  The proofs consist of three parts.
We have to derive the limit of  the expectation of the Poissonized estimator,
we have to show that the variance of the Poissonized estimator vanishes
asymptotically, and we have to prove that Poissonization is
allowed, i.e. that the difference between the original estimator
and its Poissonized version asymptotically vanishes in probability.
Here we only derive the limits of the expectation. The complete proofs are
given in Section \ref{proofs}.

We can rewrite the expectation of $ \tilde F_m(x)$ as follows
\begin{equation}
  \ex \tilde F_m(x)=\ex \frac{1}{m} \sum_{j=1}^{m}I_{[\frac{m}{n}
\bar\rho_{j,M} \leq x]}
  =
  \frac{1}{m} \sum_{j=1}^{m} P \left( \frac{m}{n}\, \bar \rho_{j,M} \leq x
  \right).
\end{equation}
Recall that for $m=M$ this gives the expectation of the Poissonized
natural estimator $\tilde F_M(x)$.

Now consider a two stage procedure.
We draw a value $z$ from the sequence of points $mq_{1,M},mq_{2,M},\dots,mq_{m,M}$ with
equal probability $1/m.$
The corresponding random variable is denoted by $Z_m$.
Note that it has distribution function $F_m$. Given $Z_m=z$ the random variable $Y_m$
is equal to $m/n$ times a Poisson$(\frac{n}{m} z)$ distributed random variable. Then
we have by conditioning on $Z_m$
\begin{equation}
 \ex \tilde F_m(x)=
\frac{1}{m} \sum_{j=1}^{m} P(\frac{m}{n}\,\bar \rho_{j,M} \leq x)
=\ex(P(Y_m \leq x|Z_m))=P(Y_m \leq x).
\end{equation}
Hence $\ex \tilde F_m(x)$ equals the distribution function of
$Y_m$ at $x$. We derive weak convergence of this distribution
function by the continuity theorem for characteristic functions.
The characteristic function of $Y_m$, denoted by $\phi_m$, is given by

\begin{equation}\label{charfu}
\phi_m(t) = \ex(e^{itY_m})=
\ex (\ex(e^{itY_m} \big| Z_m))
=\int e^{\frac{n}{m}z\left(e^{it\frac{m}{n}}-1\right)}
dF_m(z),
\end{equation}
since the characteristic function of a Poisson($\mu$) distribution
is equal to $e^{\mu\left(e^{it}-1\right)}$.
In the case of the natural estimator we have $m=M$ and hence by
(\ref{cond})
\begin{equation}\label{charone}
\phi_m(t)\to \int  e^{\lambda z\left(e^{it/\lambda}-1\right)}
dF(z),
\end{equation}
the characteristic function of the limit distribution function in (\ref{lim1}).
For the estimator based on the grouped counts we have $m/n\to 0$
by (\ref{condm}) and hence
\begin{equation}\label{chartwo}
\phi_m(t)\to \int  e^{itz}
dF(z),
\end{equation}
the characteristic function of $F$.
By the continuity theorem (\ref{charone}) and (\ref{chartwo}) imply the
conclusions of the two theorems.

\begin{rem}\label{ordered}
\textup{In Theorem \ref{thm:1.2} we can replace the condition
$F_m\stackrel{w}{\to} F$  by $F_M\stackrel{w}{\to} F$ if the
 $p_{j,M}$'s, $j=1,2,\ldots,M$ are
ordered.  A proof can be found in
Section \ref{techproofs}.  }
\end{rem}

\begin{rem}
{\rm
The condition of the weak convergence of $F_m$ to $F$ is
implied by a stronger condition in Klaassen and
Mnatsakanov (2000). Define $f_M$ by
\begin{equation}
f_M(t)=\sum_{j=1}^{M}Mp_{j,M}I_{[\frac{j-1}{M}<t\leq \frac{j}{M}]},
\hspace{3pt} 0<t\leq 1.
\end{equation}
Note that the structural distribution function $F_M$ is
the distribution function of $f_M(U)$, where $U$ is uniformly distributed
on the interval $(0,1]$. Assume that $f_M$ converges uniformly on $(0,1]$ to a
density function $f$, i.e.
\begin{equation}\label{f}
\sup_{0<t\leq 1} |f_M(t)-f(t)| \to 0.
\end{equation}
Klaassen and Mnatsakanov   proved, without requiring equal cell sizes, that this condition
implies weak consistency. Moreover, the condition (\ref{condm}) is slightly stronger then the
corresponding one   required by Klaassen and Mnatsakanov. }
\end{rem}

\bigskip

Let us consider the rate of convergence and the choice of the number of groups $m$.
Define the {\em Mean Squared Error} (MSE) of $\hat F_m(x)$ as
\begin{equation}
\mse (\hat F_m(x)) = \ex (\hat F_m(x)-F(x))^2.
\end{equation}
A standard computation shows that the mean squared error is equal to
the sum of the squared bias and
the variance.

Consider the situation where the $p_{j,M}$'s are generated by a distribution function $G$,
via
\begin{equation}\label{probs}
p_{j,M}=G(j/M)-G((j-1)/M),\quad j=1,\dots,M.
\end{equation}
Then we also have
$q_{j,M}=G(j/m)-G((j-1)/m), j=1,\dots,m$. If $G$ has a density
$g$ that is continuous and bounded then we have
\begin{equation}
m q_{j,M}=m(G(j/m)-G((j-1)/m)) = m g(\xi_{j,M}) {1\over m}= g(\xi_{j,M}),
\end{equation}
where  $\xi_{j,M}$ is a point in the interval $((j-1)/m,j/m]$.
Assuming that $g$ is also uniformly continuous on $(0,1]$ this implies $f_m(t)\to g(t)$, uniformly on
$[0,1)$. So in this situation the limit density $f$ in (\ref{f}) is
equal to $g$.

\begin{figure}[h]
 \input{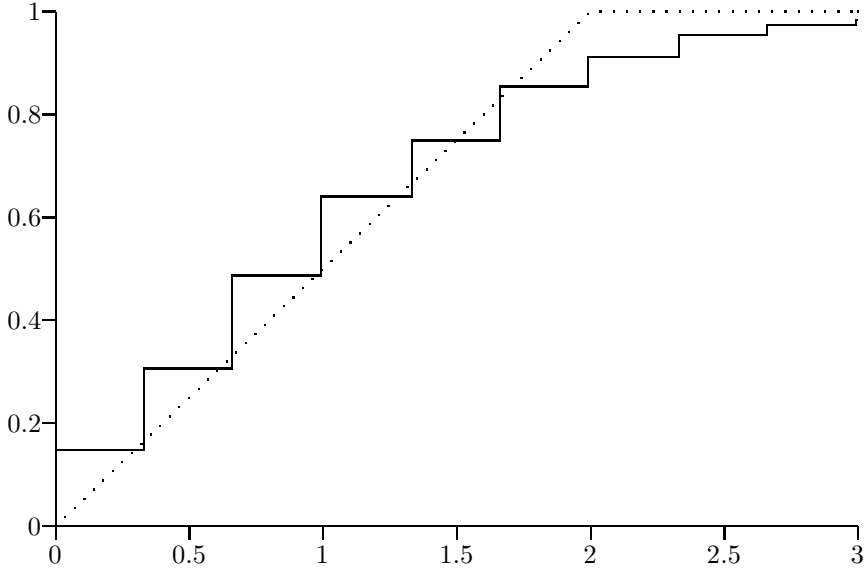}
 \caption{$\hat F_M(x)$ for $M=1000, n=3000\ (m=M=1000, k=1)$ } \label{natest}
\end{figure}

\begin{figure}[h]
 \input{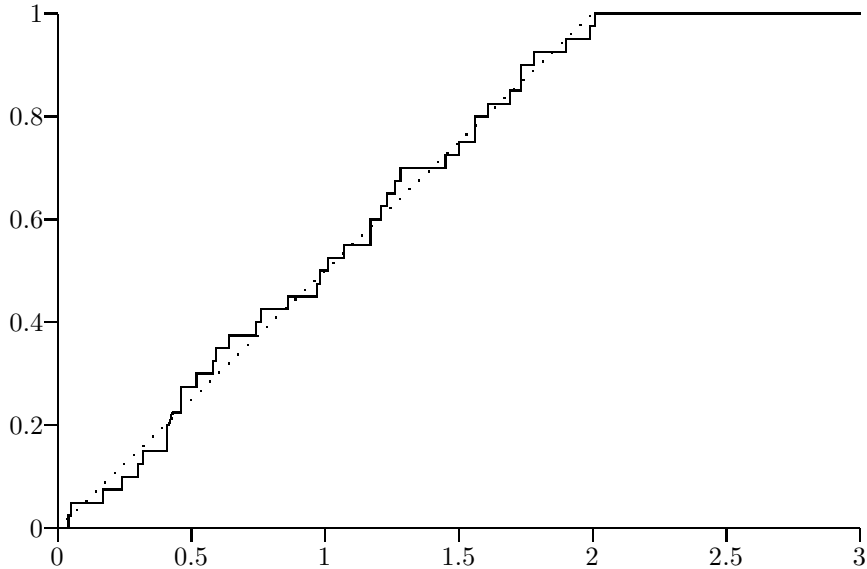}
 \caption{$\hat F_m(x)$ based on grouping with $m=40, k=25, M=1000, n=3000$} \label{k=25}
\end{figure}

\begin{figure}[h]
 \input{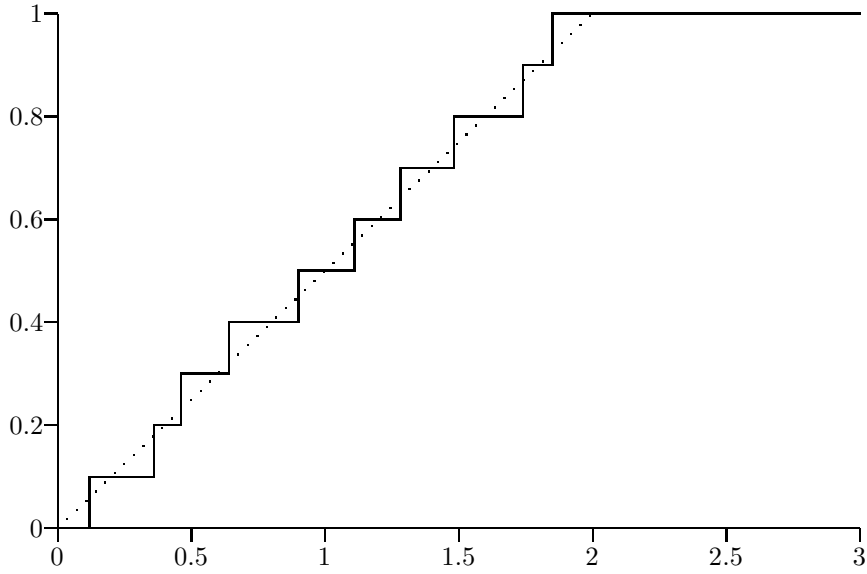}
 \caption{$\hat F_m(x)$ based on grouping with $m=10,k=100, M=1000, n=3000$} \label{k=100}
\end{figure}

Let us first present some simulation results.
Figures \ref{natest}, \ref{k=25} and  \ref{k=100} show estimates of
$F$ based on a simulated
sample where $G(x)=2x-x^2$ and $g(x)=2(1-x)$ for $0\leq x \leq 1$.
We have chosen $M=1000$  and
$n=3000$. So $\lambda$ equals three.
Since it equals the distribution function of $g(U)$,
with $U$ uniformly distributed on $[0,1)$,
the limit  structural function $F$ is given by
\begin{equation}
F(x)= \left\{
 \begin{array}{ll}
  0   & \mbox{if $x<0$}, \\
  \frac{1}{2}x & \mbox{if $0\leq x \leq 2$},\\
  1 & \mbox{if $x>2$}.
 \end{array}
\right.
\end{equation}
Figure \ref{natest} shows the result of the natural estimator.
Next we show two figures of estimates based on grouping.
In Figure \ref{k=25} we have  $k=25$ and thus $m=40$ while for Figure \ref{k=100} 
we have chosen $k=100$
and thus $m=10$.
Figure \ref{natest} shows that the natural estimator is
inconsistent, having jumps only at multiples of $1/\lambda=1/3$.
Figures \ref{k=25} and \ref{k=100} show that by grouping we achieve consistency, and
that the choice of $m$ is important. All in all the figures suggest that $k$ too 
small or too large is not wise and that there  might be an optimal cell size.

The next theorem gives some insight in  the choice of $m$. It gives
bounds on the mean squared error of $\hat F_m(x)$. These bounds
depend on $m$.

\newpage

\begin{thm}\label{thm:1.3}
Let (\ref{cond}) hold.
Assume that the cell probabilities $p_{j,M}, j=1,\dots,M$ are generated by a
distribution function $G$ as in (\ref{probs}) and that $G$ has a density that is
uniformly continuous on $(0,1]$.
Assume further that $G$ has a bounded second derivative $g$ that bounded away from zero on
$(0,1]$,
and that, for some $0<\alpha <1/6$,
\begin{equation}\label{condm2}
{n \over m (\log m)^{1/2\alpha}}\to \infty.
\end{equation}
Then we have, if $m\gg n^{1/3}$,
\begin{equation}\label{v1}
\mse (\hat F_m(x)) \leq \frac{9}{4\pi^2} (24\tau)^{4/3} \left( \frac{m}{n}\right)^{2/3}
+\frac{1}{4m}+o\left(\left( \frac{m}{n}\right)^{2/3}\right)+o\left({1\over m}\right),
\end{equation}
and if $m\ll n^{1/3}$
\begin{equation}\label{v2}
\mse (\hat F_m(x)) \leq
  \frac{1}{4m}+o\left({1\over m}\right).
\end{equation}
\end{thm}

The key idea of the proof is to exploit the fact that we have
derived the convergence of $\ex \tilde F_m(x)$, which is in fact equal
to the distribution function of $Y_m$, to $F(x)$ from the
convergence of its characteristic function $\phi_m$, cf. (\ref{charfu}), to the
characteristic function of $F$. By Esseen's smoothing lemma we get
a bound on the distance of distribution functions from the
distance of their characteristic functions. By expanding
(\ref{charfu}) we obtain a rate of convergence for the bias
$\ex \tilde F_m(x)- F(x)$ of the Poissonized estimator. The bound on the variance of the
Poissonized estimator is the same as in the proof of Theorem \ref{thm:1.2}. The
remainder of the proof consists of showing that Poissonization is
allowed in this context too.

Straightforward calculations show that the right hand side of (\ref{v1}) is
asymptotically minimized by $m_n$ if
\begin{equation}
m_n\sim  \left( \frac{\pi ^6}{6^3 (24\tau)^4} \right)^{1/5} n^{2/5}.
\end{equation}
This gives a mean squared error
\begin{equation}\label{v4}
\mse(\hat F_{m_n}(x)) \leq \frac{33}{4}\left(  \frac{(24\tau)^2}{6\pi^3}
  \right)^{2/5} n^{-2/5}+o(n^{-2/5}).
\end{equation}
The bound (\ref{v2}) of Theorem \ref{thm:1.3} gets smaller as  $m$ increases.
However, the order of
$m$ is bounded by $n^{1/3}$.
Hence, for $m\ll n^{1/3}$ we get
\begin{equation}\label{v3}
\mse (\hat F_m(x)) \gg \frac{1}{4}n^{-1/3} +o(\frac{1}{4}n^{-1/3}).
\end{equation}
Note that the bound in (\ref{v4}) is smaller than the one given in (\ref{v3}). Therefore
(\ref{v4}) gives the minimal upper bound.

\begin{rem}
{\rm
The assumption that there exists a {\em known} ordering of words in a
vocabulary, necessary for grouping,
for which (\ref{probs}) holds is not realistic. Consistent
estimators as the one in Klaassen and Mnatsakanov (2000),
which do not require such an ordering, seem to have a logarithmic
rate of convergence, as opposed to the algebraic rate in Theorem
\ref{thm:1.2}.
}
\end{rem}

\section{Proofs}\label{proofs}
\setcounter{equation}{0}

\subsection{Proof of Theorem  \ref{thm:1.1}}

The limit of $\ex \tilde F_M(x)$ is derived in the previous section. It remains to check (\ref{charone})
reformulated in the following lemma.
\begin{lem} \label{lem1}
Under the conditions of Theorem \ref{thm:1.1} we have
\begin{equation}\label{philim}
\phi_M(t)=\int  e^{\frac{n}{M}z\left(e^{it\frac{M}{n}}-1\right)} dF_M(z)
\to \int  e^{\lambda z\left(e^{it/\lambda}-1\right)}
dF(z).
\end{equation}
\end{lem}
The proof is given in Section \ref{techproofs}.

A bound on the variance of $\tilde F_M(x)$ is given by
\begin{eqnarray*}
\var (\tilde F_M(x))
&=& \var \left( \frac{1}{M} \sum_{j=1}^{M}
I_{[\frac{M}{n} \rho_{j,M} \leq x]} \right)
=
\frac{1}{M^2} \sum_{j=1}^{M}
\var \left(I_{[\frac{M}{n} \rho_{j,M} \leq x]}\right)\\
&\leq&
\frac{1}{M^2} \sum_{j=1}^{M} \frac{1}{4}
=\frac{1}{4M} \rightarrow 0.
\end{eqnarray*}
All this implies that $\tilde F_M$ is weakly consistent for
$F_{Y/\lambda}$.

Finally we show that Poissonization is allowed. We have
\begin{eqnarray*}
|\hat F_M (x) - \tilde F_M (x)|
 &=&
\Big| \frac{1}{M} \sum_{j=1}^{M}I_{[\frac{M}{n}\nu_{j,M}\leq x]}-
 \frac{1}{M} \sum_{j=1}^{M}I_{[\frac{M}{n}\rho_{j,M}\leq x]} \Big |\\
&\leq&
\frac{1}{M} |N-n| = \frac{n}{M} \Big|\frac{N}{n}-1 \Big| \rightarrow
0,
\end{eqnarray*}
almost surely and in probability.
This implies that $\hat F_M$ is weakly consistent for
$F_{Y/\lambda}$ too as stated in the theorem.

\subsection{Proof of Theorem \ref{thm:1.2}}

The limit of $\ex \tilde F_m(x)$ is derived in the previous section. It remains to check (\ref{chartwo})
reformulated in the following lemma.
\begin{lem} \label{lem2}
Under the conditions of Theorem \ref{thm:1.2} we have
\begin{equation}
\phi_m(t)=\int e^{\frac{n}{m}z\left(e^{it\frac{m}{n}}-1\right)} dF_m(z) \to
\int e^{itz}dF(z).
\end{equation}
\end{lem}
The proof can be found in Section \ref{techproofs}.

Here we bound the variance of $\tilde F_m(x)$ as follows
\begin{eqnarray}
\var \tilde F_m(x) &=& \var \frac{1}{m}
\sum_{j=1}^{m}I_{[\frac{m}{n} \bar \rho_{j,M} \leq x]}
=
\frac{1}{m^2} \sum_{j=1}^{m} \var I_{[\frac{m}{n} \bar \rho_{j,M} \leq
x]}\label{varbound}\\
&\leq&
\frac{1}{m^2} \sum_{j=1}^{m} \frac{1}{4}
=
\frac{1}{4m} \to 0.\nonumber
\end{eqnarray}
This implies that $\tilde F_m(x)$ is a weakly consistent for
$F(x)$.

In order to transfer the weak consistency result to the original estimator
we must show that  we may indeed Poissonize, i.e. we must show that
$\hat F_m (x) - \tilde F_m (x)$ vanishes in probability.

We need the Bernstein inequality for Poisson random variables.
If $X$ has a Poisson distribution then
\begin{equation}\label{poisbound}
P\Big({{|X-\ex X|}\over (\ex X)^{1/2}}\geq\epsilon\Big)
\leq
2\exp\Big(-\ {\epsilon^2\over{2+\epsilon (\ex X)^{-1/2}}}\Big),
\end{equation}
cf. Lemma 8.3.4 in Reiss (1993). It also follows from Inequality 1 on page 485 of
Shorack and Wellner (1986).

Write $z_{j,n}=m q_{j,M}$.
Note that, since the distributions induced by the $F_m$ are concentrated
on a bounded set,
we have $\max_{1\leq j\leq m} z_{j,n}\leq c$
for some constant $c>0$.
Hence, for all $\delta >0$, we have
\begin{eqnarray}
\lefteqn{\sum_{j=1}^m
P\Big(\Big|{m\over n}\,\bar\rho_{j,M}-z_{j,n}\Big|
\geq\delta\Big)= }\nonumber\\
&=&
\sum_{j=1}^m
P\Big(|\bar\rho_{j,M}-{n\over m}\ z_{j,n}|
\geq{n\over m}\ \delta\Big)\nonumber\\
&=&
\sum_{j=1}^m
P\Big({{|\bar\rho_{j,M}-nq_{j,M}|}\over
{(nq_{j,M})^{1/2}}}
\geq\Big({n\over m}\Big)^{1/2}{1\over \sqrt{z_{j,n}}}\
\delta\Big)\label{expbound}\\
&\leq&
\sum_{j=1}^m
2\exp\Big(-\delta^2\ {n\over m}\ {1\over z_{j,n}}\
{1\over {2 + \delta{1\over z_{j,M}}}}\Big)\nonumber\\
&\leq&
2 m\exp\Big(-\delta^2\ {n\over m}\
{1\over {2c + \delta}}\Big)\nonumber\\
&=&
2\exp\Big( \log m \Big(-{n\over m\log m}\
{\delta^2\over {2c + \delta}}+1\Big)\Big)\to 0,\nonumber
\end{eqnarray}
by (\ref{condm}).

By $\max_{1\leq j\leq m} q_{j,M}\to 0$ we have
$\var(\bar\nu_{j,M})=nq_{j,M}(1-q_{j,M})\sim
nq_{j,M}=\var(\bar\rho_{j,M})$. By the Bernstein inequality for
binomial random variables, cf. Shorack and Wellner (1986), p 440, it now follows that
for $\delta>0$
\begin{equation}\sum_{j=1}^m
P\Big(\Big|{m\over n}\,\bar\nu_{j,M}-z_{j,n}\Big|
\geq\delta\Big)\to 0.
\end{equation}
This implies that with probability approaching one we have
\begin{equation}\label{bounds}
\Big|{m\over n}\,\bar\nu_{j,M}-z_{j,n}\Big|
<\delta
\quad\mbox{and}\quad
\Big|{m\over n}\,\bar\rho_{j,M}-z_{j,n}\Big|
<\delta,\quad j=1,\dots,m.
\end{equation}
Consequently, (\ref{bounds}) implies
\begin{equation}\label{enclosed}
{1\over m}\sum_{j=1}^m
(\I_{[z_{j,n}\leq x-\delta]} - \I_{[z_{j,n}\leq x+\delta]})
\leq
\hat F_m (x) - \tilde F_m (x)
\leq
{1\over m}\sum_{j=1}^m
(\I_{[z_{j,n}\leq x+\delta]} - \I_{[z_{j,n}\leq x-\delta]}).
\end{equation}
By the weak convergence of $F_m$ to $F$, if $x-\delta$ and $x+\delta$ are continuity points of $F$,
 the left and right hand side converge to
$F(x-\delta)-F(x+\delta)$ and $F(x+\delta)-F(x-\delta)$
respectively. Now, for given $\epsilon>0$, choose $\delta$
such that $F(x+\delta)-F(x-\delta)$
is smaller than $\epsilon$ and we have shown
\begin{equation}
P(|\hat F_m (x) - \tilde F_m (x)|\geq \epsilon)\to 0.
\end{equation}
Hence $\hat F_m (x) - \tilde F_m (x)$ vanishes in probability, proving
that Poissonization is allowed.

\subsection{Proof of Theorem \ref{thm:1.3}}

First we consider the mean squared error of the Poissonized estimator.
By a standard calculation  we have
\begin{equation}
\mse (\tilde F_m(x)) = (\ex \tilde F_m(x)-F(x))^2 + \var(\tilde F_m(x))
\end{equation}
A bound on the
variance is already given by (\ref{varbound}).
It is harder to obtain a bound on the bias. We shall use the convergence of
the
characteristic function of $Y_m$ to the characteristic function of $Y$ in the
proof of Theorem~\ref{thm:1.1} and
Esseen's smoothing
lemma, see Feller (1966), Section XIV 3, Lemma 2 on page 538.

\begin{lem}[Esseen's smoothing lemma]
Let $F$ be a probability distribution function with vanishing expectation
and characteristic
function $\varphi$. Suppose $F-G$ vanishes at $\pm \infty$ and that $G$ has
a derivative $g$
such that $|g|\leq\tau$. Finally, suppose that $g$ has a continuously
differentiable Fourier
transform $\gamma$ such that $\gamma(0)=1$ and $\gamma'(0)=0$. Then, for
all $x$ and $T>0$
\begin{equation}
|F(x)-G(x)|\leq {1\over\pi} \int_{-T}^T
\Big|{{\varphi(t)-\gamma(t)}\over t}\Big|dt
+ {24 \tau\over \pi T}.
\end{equation}
\end{lem}

Now apply this lemma with $F$ equal to the distribution function of $Y_m$
and $G$ equal
to the limit structural distribution function $F$.
Note that both distribution functions have expectation one and that the induced
distributions are concentrated on $[0,\infty)$. Then
\begin{equation} \label{esl}
|\ex \tilde F_m(x) - F(x)|
\leq
{1\over\pi} \int_{-T}^T
\Big|{{ \ex e^{itY_m}-\ex e^{itZ}}\over t}\Big|dt
+ {24\tau\over \pi  T}.
\end{equation}
Let us first consider the integrand. Write
\begin{eqnarray}
\lefteqn{|\ex e^{itY_m}-\ex e^{itZ}|=
|\int e^{\frac{n}{m}z\left(e^{it\frac{m}{n}}-1\right)} dF_m(z) -  \int
e^{itz}dF(z)|}
\nonumber\\
&\leq&
|\int e^{\frac{n}{m}z\left(e^{it\frac{m}{n}}-1\right)} dF_m(z)-
\int e^{itz} dF_m(z)|\label{w1}\\
&+&|\int e^{itz} dF_m(z)-  \int e^{itz}dF(z)|
\label{w2}
\end{eqnarray}
For $n$ large we have
\begin{equation}
e^{\frac{n}{m}z(e^{it\frac{m}{n}}-1)}=e^{\frac{n}{m}z(1+it\frac{m}{n}-\frac{1}{2
}t^2(\frac{m}{n})^2+R_n(t)-1)}
=e^{itz-\frac{1}{2}\frac{m}{n}t^2z+R_n(t)\frac{n}{m}z},
\end{equation}
where $R_n(t)=e^{it\frac{m}{n}}-1-it\frac{m}{n}+\frac{1}{2}t^2\frac{m^2}{n^2}$.
Note that
\begin{eqnarray}\label{Rn}
|R_n(t)|\leq \frac{1}{6}t^3\frac{m^3}{n^3}
\end{eqnarray}
and that  for $w \epsilon {\mathbb C}$, and $|w|$ small enough, we have
\begin{equation}
|e^w-1| \leq 4|w|.
\end{equation}
Hence
\begin{equation}\label{w3}
\big|e^{\frac{n}{m}zR_n(t)}-1 \big| \leq \frac{2}{3}|z||t^3|\frac{m^2}{n^2}.
\end{equation}
So we can bound the term (\ref{w1}) as follows
\begin{eqnarray} \label{w4}
\lefteqn{\Big|\int e^{\frac{n}{m}z\left(e^{it\frac{m}{n}}-1\right)} dF_m(z)-
\int e^{itz} dF_m(z)\Big|} \nonumber\\
&\leq&
\Big|\int \left( e^{itz-\frac{1}{2}\frac{m}{n}t^2z}-e^{itz} \right)dF_m(z) \Big|
 \nonumber\\
&+&
\Big|\int \left( e^{\frac{n}{m}z\left(e^{it\frac{m}{n}}-1\right)}-e^{itz-\frac{1}{2}
\frac{m}{n}t^2z}\right)dF_m(z)\Big|\\
&\leq&
\int \Big|e^{-\frac{1}{2}\frac{m}{n}t^2z} -1\Big|dF_m(z) +\int \Big|e^{\frac{n}{m}
zR_n(t)} -1\Big|dF_m(z) \nonumber\\
&\leq&
\frac{1}{2}{m\over n}t^2\int zdF_m(z) + \frac{2}{3}\frac{m^2}{n^2}|t|^3 \int zdF_m(z).
 \nonumber\\
&=&
\frac{1}{2}{m\over n}t^2 + \frac{2}{3}\frac{m^2}{n^2}|t|^3.
 \nonumber
\end{eqnarray}
For (\ref{w3}) to hold we have tacitly assumed that $(n/m)zR_n(t)$ vanishes for
$-T\leq t \geq T$. By (\ref{Rn})
and the fact that  $Z_m$ is almost surely bounded by the same
 constant for all $m$, it suffices
to check that $(m^2/n^2)t^3 \to 0$ for $-T\leq t \geq T$. Further on in the proof $T$
will depend on $n$. The condition
is satisfied for our two choices of $T_n$ in (\ref{Tn1}) and (\ref{Tn2}).

For the first term in (\ref{esl}) we get
\begin{eqnarray} \label{w5}
\lefteqn{\frac{1}{\pi}\int_{-T}^T\frac{1}{|t|}\Big|\int e^{\frac{n}{m}z
\left(e^{it\frac{m}{n}}-1\right)} dF_m(z)-\int e^{itz} dF_m(z)\Big|dt} \nonumber\\
&\leq&
{1\over 2\pi} {m\over n} \int_{-T}^T|t|dt+{2\over 3\pi} {m^2\over n^2} \int_{-T}^T t^2
 dt\\
&=&
{1\over 2\pi} {m\over n} T^2+{4\over 9\pi} {m^2\over n^2} T^3 .\nonumber
\end{eqnarray}

Let the function $f_m$ be defined by
\begin{equation}
f_m(t)=\sum_{j=1}^{m}mq_{j,m}I_{[\frac{j-1}{m}<t\leq \frac{j}{m}]},
\hspace{3pt} 0<t\leq 1.
\end{equation}
Then  $F_m$ is the distribution function of $f_m(U)$ where $U$ is uniformly
distributed on $(0,1]$. Since $f_m$ converges uniformly to $g$ the limit
distribution function $F$ is the distribution function of $g(U)$.
Hence
\begin{equation}
\int e^{itz} dF_m(z)-  \int e^
{itz}dF(z)=\int_0^1( e^{itf_m(u)}-e^{itg(u)})du.
\end{equation}
Integrated over the intervals
$((j-1)/m,j/m]$, the constant $mq_{j,M}$ yields the same value as
$g$ integrated over these intervals. So we can write
\begin{eqnarray*}
\lefteqn{  e^{itf_m(u)}-e^{itg(u)}= e^{itf_m(u)} \left( 1-e^{it(g(u)-f_m
(u))}  \right)}\\
&=&
e^{itf_m(u)} \left(  it(f_m(u)-g(u))+R_n(t,u)
\right),
\end{eqnarray*}
where
\begin{equation}
|R_n(t,u)|\leq \frac{1}{2} t^2 (g(u)-f_m(u))^2
\end{equation}
And hence, if $g$ has a bounded derivative on $(0,1]$,
\begin{eqnarray*}
\lefteqn{\Big|\int \left(  e^{itf_m(u)}-e^{itg(u)} \right) du\Big|
=
\Big| \int e^{itf_m(u)} R_n(t,u) du\Big|}\\
&\leq&
\frac{1}{2} t^2 \int (f_m(u)-g(u))^2 du \leq \frac{c}{2} \frac{t^2}{m^2},
\end{eqnarray*}
where $c$ is a positive constant.
This implies
\begin{equation} \label{z1}
\frac{1}{\pi}\int_{-T}^T\frac{1}{|t|}\Big|\int e^{itz} dF_m(z)-  \int e^
{itz}dF(z)\Big|dt
\leq
\frac{1}{\pi}\int_{-T}^T\frac{1}{|t|}\frac{c}{2} \frac{t^2}{m^2}\,dt =\frac{c}{2
\pi}\frac{T^2}{m^2}.
\end{equation}
Hence, for all $x$ and $T>0$
\begin{equation} \label{z2}
|\ex \tilde F_m(x) - F(x)|
\leq
\frac{4}{9\pi}\frac{m^2}{n^2}T^3+{1\over 2\pi} {m\over n} T^2+\frac{c}{2\pi}\frac{1}{m^2}T^2+\frac{24\tau}{\pi T}.
\end{equation}

First assume  that $m\gg n^{1/3}$. Then equation (\ref{z2}) becomes asymptotically
\begin{equation} \label{z3}
|\ex \tilde F_m(x) - F(x)|
\leq
{1\over 2\pi} {m\over n} T^2+\frac{24\tau}{\pi T}.
\end{equation}
The value $T_n$ that minimizes the right hand side of (\ref{z3})
is given by
\begin{equation}\label{Tn1}
T_n=(24 \tau)^{1/3}\left( \frac{n}{m}  \right) ^{1/3}.
\end{equation}
Hence the bias can be asymptotically bounded by
\begin{equation}
|\ex \tilde F_m(x) - F(x)| \leq \frac{3}{2\pi} (24\tau)^{2/3} \left( \frac{m}{n}\right)^{1/3}
+o\Big(\left( \frac{m}{n}\right)^{1/3}\Big)
\end{equation}
and the mean squared error by
\begin{equation}\label{z4}
\mse (\tilde F_m(x)) \leq \frac{9}{4\pi^2} (24\tau)^{4/3} \left( \frac{m}{n}\right)^{2/3} +\frac{1}{4m}
+o\Big(\left( \frac{m}{n}\right)^{2/3}\Big)+o\Big({1\over m}\Big).
\end{equation}

If $m\ll n^{1/3}$, by minimizing the third and fourth term in
(\ref{z2}),
we get, by choosing
\begin{equation}\label{Tn2}
T_n=c^{-1/3}(24 \tau)^{1/3}m^{2/3},
\end{equation}
that asymptotically
\begin{equation}\label{z5}
\mse (\tilde F_m(x)) \leq \frac{9}{4\pi^2}\,c^{2/3} (24\tau)^{4/3} m^{-4/3} +\frac{1}{4m}
+ o( m^{-4/3})+o\Big({1\over m}\Big)={1\over 4m} +o\Big({1\over m}\Big).
\end{equation}

We have now derived the asymptotic bounds on the mean squared
error of the Poissonized estimator. We will show that Poissonization is
allowed. By the triangle inequality we have
\begin{equation}\label{mseest}
\mse(\hat F_m(x))^{1/2}\leq \mse(\tilde F_m(x))^{1/2} +
(\ex(\hat F_m(x) -\tilde F_m(x))^2)^{1/2}.
\end{equation}
The second term on the right hand side can be dealt with using
the following lemma. Its proof is given in Section
\ref{techproofs}.

\begin{lem}\label{msepois}
Under the conditions of Theorem \ref{thm:1.3} and we have
for any $0<\alpha<{1\over 6}$
\begin{equation}\label{p2}
\ex ({\tilde F}_m(x) -   {\hat F}_m(x))^2 =
O\Big(\Big({m\over n}\Big)^{1  - 2\alpha}\Big)
+O\Big({1\over m^2}\Big).
\end{equation}
\end{lem}
By this order bound and (\ref{z4}) and (\ref{z5}) it follows 
that $(\ex(\hat F_m(x) -\tilde F_m(x))^2)^{1/2}$ is
asymptotically negligible compared to $\mse(\tilde F_m(x))^{1/2}$.
Hence Poissonization is allowed.

\section{Technical proofs} \label{techproofs}
\setcounter{equation}{0}
\subsection{Proof of Lemma \ref{lem1}}
\noindent
Recall that $(n/M)Y_M$, given $Z_M=z$, has a Poisson($\frac{n}{M} z$)
distribution. We have
$Z_M
\stackrel{w}{\to} Z$, so  $F_M(w) \to F(w)$ at all continuity points $w$  of
$F.$ Let $\psi_M$ denote the characteristic function   of $(n/M)Y_M$.
Then
\begin{equation}
\psi_M(t) = \ex \left( e^{it(n/M)Y_M} \right)
=
\ex \left( \ex \left( e^{it(n/M)Y_M}|Z_M \right) \right)
=
\int_{-\infty}^{\infty} e^{\frac{n}{M}z \left(e^{it}-1 \right)} dF_M(z).
\end{equation}
Consider t fixed. For $z \in [0,w]$ we have
$$
\big| e^{\frac{n}{M}z \left(e^{it}-1 \right) }-
e^{\lambda z \left( e^{it}-1 \right)} \big|
=
\big| e^{\frac{n}{M} z \left( e^{it}-1\right)}\big| \,
\Big| 1-e^{\left( \lambda -\frac{n}{M} \right) z \left( e^{it}-1
\right)} \Big|
\leq
 \big| 1-e^{\left( \lambda -\frac{n}{M} \right) z \left( e^{it}-1
\right)} \big| \to 0
$$
or equivalently
\begin{eqnarray*}
  e^{\frac{n}{M}z\left(e^{it}-1\right)} \to e^{\lambda
z\left(e^{it}-1\right)}.
\end{eqnarray*}
This also holds for $z$ replaced by
$z_n$, for every sequence \{$z_n$\} with values in $[0,w]$, showing that
the convergence is uniform in $z$.
Hence for $\epsilon >0$ and $n$ large enough
\begin{equation}  \label{x1}
  \big| e^{\frac{n}{M}z\left(e^{it}-1\right)} -
        e^{\lambda z\left(e^{it}-1\right)}      \big|
  \leq \epsilon/2
\end{equation}
for all $z\in [0,w].$

Let $w$ be a continuity point of $F$.
Note that, because $F$ and $F_M$ vanish on the negative half line, the point -1
is also a continuity point.
Then, according to the Helly-Bray theorem and because characteristic
functions are continuous, we can conclude that
\begin{equation}
  \int_{-1}^{w} e^{\lambda z \left( e^{it}-1 \right)} dF_M(z) \to
  \int_{-1}^{w} e^{\lambda z \left( e^{it}-1 \right)} dF(z).
\end{equation}
So for $n$ large enough \\
\begin{equation}  \label{x2}
  \Big|   \int_{-1}^{w} e^{\lambda z \left( e^{it}-1 \right)} dF_M(z) -
  \int_{-1}^{w} e^{\lambda z \left( e^{it}-1 \right)} dF(z) \Big|
  \leq \epsilon/2.
\end{equation}
Because of (\ref{x1}) and (\ref{x2}) we now  have

\begin{eqnarray}
  \lefteqn{ \Big|   \int_{-1}^{w} e^{\frac{n}{M} z \left( e^{it}-1 \right)}
dF_M(z)
    -
    \int_{-1}^{w} e^{\lambda z \left( e^{it}-1 \right)} dF(z) \Big| }\nonumber\\
     &\leq&
    \Big| \int_{-1}^{w} \left(e^{\frac{n}{M} z \left(e^{it}-1 \right)}
    -e^{\lambda z \left(e^{it}-1 \right)} \right) dF_M(z) \Big|   \label{y1}\\
      &+&
    \Big| \int_{-1}^{w} e^{\lambda z \left(e^{it}-1 \right)} dF_M(z)
    -\int_{-1}^{w} e^{\lambda z \left(e^{it}-1 \right)}dF(z) \Big|\nonumber\\
     &\leq& \epsilon.\nonumber
\end{eqnarray}

Next choose the continuity point $w$ such that
$1-F(w)<\epsilon/2$.
Since $F_M(w) \to F(w)$ we also have $0\leq 1-F_M(w)<\epsilon/2$, for $n$ large enough.
This implies
\begin{eqnarray}
\lefteqn{\Big| \int_{w}^{\infty} e^{\frac{n}{M}z \left(e^{it}-1\right)}dF_M(z)-
\int_{w}^{\infty} e^{\lambda z \left(e^{it}-1\right)}dF(z)\Big|}\nonumber  \\
&\leq&
\int_{w}^{\infty} \big| e^{\frac{n}{M}z \left(e^{it}-1\right)} \big|
dF_M(z)+
\int_{w}^{\infty} \big| e^{\lambda z \left(e^{it}-1\right)} \big| dF(z)
\label{y2} \\
&=&
1-F_M(w)+1-F(w)<\epsilon.\nonumber
\end{eqnarray}
The inequalities (\ref{y1}) and (\ref{y2})
show that
\begin{equation}
\Big|
\int^{\infty}_{-\infty}e^{\frac{n}{M}z\left(e^{it}-1\right)}dF_M(z) -
\int^{\infty}_{-\infty}e^{\lambda z\left(e^{it}-1\right)}dF(z) \Big|
<2\epsilon
\end{equation}
for $n$ large enough. Hence
\begin{equation}
\int^{\infty}_{-\infty}e^{\frac{n}{M}z\left(e^{it}-1\right)}dF_M(z) \to
\int^{\infty}_{-\infty}e^{\lambda z\left(e^{it}-1\right)}dF(z) .
\end{equation}
Since convergence of characteristic functions is uniform on bounded
intervals we also have (\ref{philim}).

\subsection{Proof of Lemma \ref{lem2}}
\noindent
The proof is similar to the proof in the previous section.
Note that
\begin{equation}\label{b2}
  \lim_{n \to \infty}  \frac{n}{m} z  \left( e^{it\frac{m}{n}}-1  \right)
 =itz ,
\end{equation}
 uniformly for $z \epsilon[-1,w]$.
Let $\epsilon>0$  and $w$  be  a continuity point of $F$.
By the Helly-Bray
theorem and (\ref{b2}), we have, for $n$ large enough,
\begin{eqnarray}
  \lefteqn{\big| \int_{-1}^{w}
e^{\frac{n}{m}z\left(e^{it\frac{m}{n}}-1\right)}dF_m(z)
  -
  \int_{-1}^{w} e^{itz}dF(z)  \big|} \nonumber  \\
  &\leq&
  \big| \int_{-1}^{w}
e^{\frac{n}{m}z\left(e^{it\frac{m}{n}}-1\right)}dF_m(z)
  -
  \int_{-1}^{w} e^{itz}dF_m(z)
\big|\label{b3}  \\
  &+&
  \big| \int_{-1}^{w}
e^{itz}dF_m(z)
  -
  \int_{-1}^{w} e^{itz}dF(z)  \big| \nonumber \\
  &<& \epsilon.\nonumber
\end{eqnarray}
Now choose  $w$  such that   $1-F(w)<\epsilon/2$. Then, for $n$ large enough,
\begin{eqnarray}
  \lefteqn{\left|\int_{w}^{\infty}
e^{\frac{n}{m}z\left(e^{it\frac{m}{n}}-1\right)}dF_m(z)
  -
  \int_{w}^{\infty} e^{itz}dF(z)  \right|}\nonumber   \\
  &\leq&
  \int_{w}^{\infty} \big|e^{\frac{n}{m}z\left(e^{it\frac{m}{n}}-1\right)}
\big| dF_m(z)
  +
  \int_{w}^{\infty} |e^{itz}|dF(z) \label{b1}   \\
  &\leq&
  1-F_m(w)+1-F(w)<\epsilon.\nonumber
\end{eqnarray}
As in the previous section the inequalities (\ref{b1}) and (\ref{b3}) prove
the lemma.

\subsection{Proof of Remark \ref{ordered}}
\noindent
We assume that the set of the $p_{j,M}$'s is ordered. So
$p_{1,M} \leq p_{2,M} \leq \cdots \leq p_{M,M}$.
Let $x$ be a continuity point of $F$. We want to show that
\begin{eqnarray} \label{diff}
   \big|F_M(x)-F_m(x)\big|=
   \big|\frac{1}{M}\sum_{j=1}^{M}I_{[Mp_{j,M}\leq
x]}-\frac{1}{m}\sum_{i=1}^{m}I_{[mq_{i,M}\leq x]} \big|
\end{eqnarray}
vanishes since this implies  that $F_m(x) \to F(x)$
follows from $F_M(x)\to F(x)$.

Assume that in the $\beta$ first groups of the $m$, $\beta =0,\dots,m$,
we have $Mp_{j,M}\leq x$ and that in the ($\beta$+1)th
group for the first $\alpha$, $\alpha=1,\dots,k$, of the $p_{j,M}$'s we have $Mp_{j,M}\leq x$
and that for the others $Mp_{j,M}> x$.
Then in total exactly $k\beta+\alpha$ of the $p_{j,M}$'s satisfy
$Mp_{j,M}\leq x$. Note that both $\beta$ and $\alpha$ depend on $M$ and
$x$.

Let us focus   on the $i$-th group, where $i=0,\dots,  \beta$. Then we have
$Mp_{j,M} \leq x$
for all $j=k_{i-1}+1,\dots,k_{i-1}+k=k_i$ and
hence for all $i=1,\dots,\beta$
\begin{equation}
   M \sum_{j=k_{i-1}+1}^{k_i}p_{j,M} \leq kx.
\end{equation}
This implies $mq_{i,M} \leq x$.

We can now bound the
difference (\ref{diff}).  We get
\begin{eqnarray*}
\lefteqn{\big|\frac{1}{M}\sum_{j=1}^{M}I_{[Mp_{j,M}\leq
x]}-\frac{1}{m}\sum_{i=1}^{m}I_{[mq_{i,M}\leq x]} \big|}\\
&=& \left| \frac{k\beta+\alpha}{M}-\frac{\beta }{m}-
\frac{1}{m} I_{[mq_{\beta+1,M}\leq x]} \right|
=\left| \frac{\alpha}{M}-\frac{c}{m} \right| \to 0,
\end{eqnarray*}
since  $\alpha/M \leq k/M \to 0$.

\subsection{Proof of Lemma \ref{msepois}}

Let $\delta_n= (m/n)^{1/2-\alpha}$ and let $A_n$ denote the event
\begin{equation}\label{boundsn}
\Big|{m\over n}\,\bar\nu_{j,M}-z_{j,n}\Big|
<\delta_n
\quad\mbox{and}\quad
\Big|{m\over n}\,\bar\rho_{j,M}-z_{j,n}\Big|
<\delta_n,\quad j=1,\dots,m.
\end{equation}
Then, as in   (\ref{expbound}) we have, for $n$ large enough
\begin{eqnarray*}
P(A_n^c)&\leq&\sum_{j=1}^n\Big\{
P\Big(\Big|{m\over n}\,\bar\nu_{j,M}-z_{j,n}\Big|\geq
\delta_n\Big)
+
P\Big(\Big|{m\over n}\,\bar\rho_{j,M}-z_{j,n}\Big|\geq\delta_n\Big)\Big\}\\
&\leq&
4 m \exp\Big(-\delta_n^2\, {n\over m} \, {1\over
2c+\delta_n}\Big)\\
&=&
4  \exp\Big(-\Big({n\over m}\Big)^{2\alpha} \, {1\over
c} -\log m\Big)\\
&\leq&
4  \exp\Big(-\log m \Big( {1\over
c}{n^{2\alpha}\over {m^{2\alpha}\log m}} -1\Big)\\
&\leq&
{1\over m^2}.
\end{eqnarray*}
Using (\ref{enclosed}) we write
\begin{eqnarray*}
\lefteqn{\ex ({\tilde F}_m(x) -  {\hat F}_m(x))^2}\\
&=&
\ex ({\tilde F}_m(x) -  {\hat F}_m(x))^2\I_{A_n} +
\ex ({\tilde F}_m(x) -  {\hat F}_m(x))^2\I_{A_n^c}\\
&\leq&
\Big({1\over m}\sum_{j=1}^m
(\I_{[z_{j,n}\leq x+\delta_n]} - \I_{[z_{j,n}\leq x-\delta_n]})\Big)^2
+P(A_n^c)\\
&=&
(F_m(x+\delta_n) - F_m(x-\delta_n))^2+P(A_n^c).
\end{eqnarray*}
Now recall that $F_m$ is the empirical distribution function based
on the values $mq_{j,M}, j=1,\dots,m$. If $g'(x)>0$ then each of
these values are order $1/m$ apart. Hence there are order
$\delta_n/(1/m)=m\delta_n$ values in the interval $(x-\delta_n,
x+\delta_n]$, each contributing $1/m$ to the probability.
So
\begin{equation}\label{probbound}
F_m(x+\delta_n) - F_m(x-\delta_n)=O(\delta_n).
\end{equation}
Hence
\begin{equation}
\ex ({\tilde F}_m(x) -  {\hat F}_m(x))^2
=
O(\delta_n^2)+O({1\over m^2}),
\end{equation}
which completes the proof of the lemma.

\vspace{2 cm}

\noindent{\bf\Large References}

\begin{description}

\item Feller, W. (1966), {\it An Introduction to Probability Theory and Its
Applications},
Wiley, New York.

\item
Klaassen, C.A.J. and R.M. Mnatsakanov (2000) {\it Consistent estimation
of the structural distribution function\/},  Scand. J. Statist.,
27, 733--746.

\item
Reiss, R.-D. (1993) {\it A Course on Point Processes}, Springer-verlag,
New York.

\item
Shorack, G.R. and J.A. Wellner (1986) {\it Empirical Processes with Applications to
Statistics\/},
Wiley, New York.

\end{description}
\end{document}